\documentclass[10pt]{amsart}
\usepackage{url,graphicx}
\usepackage[dvips]{epsfig}
\usepackage{latexsym,color}
\usepackage{amscd,amssymb,verbatim}
\newcommand{\be}{\begin{enumerate}}
\newcommand{\ee}{\end{enumerate}}

\newcommand{\thom}{\text{\tt Hom}\,}
\newcommand{\ra}{\rightarrow}
\newcommand{\nin}{\noindent}
\newcommand{\pr}{\noindent{\bf Proof. }}

\newcommand{\sv}{C_{2r+1}}
\newcommand{\zz}{{\mathbb Z}_2}
\newcommand{\cob}{\partial}
\newcommand{\hck}{X_{r,n}}

\newcommand{\rp}{{\mathbb R\mathbb P}}

\newtheorem{thm}{Theorem}[section]
\newtheorem{df}[thm]{Definition}
\newtheorem{lm}[thm]{Lemma}

\numberwithin{equation}{section}
\numberwithin{figure}{section}

\begin{document}

\title
{Cobounding odd cycle colorings}

\author{Dmitry N. Kozlov}
\address{Department of Computer Science, ETH Z\"urich, Switzerland}
\email{dkozlov@inf.ethz.ch}
\thanks {Research supported by Swiss National Science Foundation Grant PP002-102738/1}


\begin{abstract}
We give a~very short self-contained combinatorial proof of the Babson-Kozlov conjecture, by presenting a~cochain whose coboundary is the desired power of the characteristic class.
\end{abstract}

\maketitle

\section{Preliminaries}

The study of the following family of complexes has recently been undertaken in connection with equivariant obstructions to graph colorings.

\begin{df} \label{df:hom}
For any graphs $T$ and $G$, $\thom(T,G)\subseteq\prod_{x\in V(T)} \Delta^{V(G)}$ consists of all cells $\sigma=\prod_{x\in V(T)}\sigma_x$, such that for any $x,y\in V(T)$, if $(x,y)\in E(T)$, then $(\sigma_x,\sigma_y)$ is a~complete bipartite subgraph of~$G$.
\end{df}
\nin In particular, the cells of $\thom(T,G)$ are indexed by functions $\sigma:V(T)\ra 2^{V(G)}$ satisfying that additional property, and $\dim\sigma=\sum_{v\in V(T)}(|\sigma(v)|-1)$. We refer the reader to the survey \cite{IAS} for an~introduction to the subject of $\thom$ complexes.

The study of the complexes $X_{r,n}:=\thom(C_{2r+1},K_n)$, $n\geq 3$, has been of special interest. Here for $r\in{\mathbb N}$, we let $\sv$ denote both the cyclic graph with $2r+1$ vertices and the additive cyclic group with $2r+1$ elements. The adjacent vertices of $v\in\sv$ get labels $v+1$ and $v-1$. Taking the negative in the cyclic group gives an involution $\gamma$ of the graph with a~fixed vertex $0$ and a~flipped edge $(r,r+1)$. Then $(X_{r,n},\gamma)$ is a~$\zz$-space, hence the Stiefel-Whitney characteristic class $w_1(X_{r,n})\in H^1(X_{r,n}/\zz;\zz)$ of the associated line bundle can be considered.

\begin{thm} \label{thm:main} {\rm (Babson-Kozlov conjecture).}
\nin We have $w_1^{n-2}(X_{r,n})=0$.
\end{thm}

The case $r=1$ was settled in \cite{BK03b}. For $r\geq 2$, and odd $n$, it was proved in \cite{BK03c}, see also \cite{IAS}, where the remaining case: $r\geq 2$, $n\geq 4$, $n$ is even, was conjectured. The latter was then proved in \cite{Sch,Sch3}. In the next section we give a~short self-contained combinatorial proof of Theorem~\ref{thm:main} covering all cases: we simply take a~cochain representative of $w_1^{n-2}(X_{r,n})$ and certify that it is a~coboundary.

First we fix notations. For $t\in{\mathbb N}$, we set $[t]:=\{1,\dots,t\}$. For a cell complex $X$, we let $X^d$ denote the set of $d$-dimensional cells of $X$. Since we are working over $\zz$, we may identify $d$-cochains with their support subsets of $X^d$. Then, the cochain addition is replaced by the symmetric difference of sets, denoted $\oplus$. For $S\subseteq\hck^d$ the coboundary operator translates to $\cob S=\oplus_{\sigma\in S}\{\tau\in\hck^{d+1}\,|\,\tau\supset\sigma\}$.

 If $r$ is even, set $t:=r/2$, and $v_i:=r-2i+1$, for $i\in[t]$, else set $t:=(r+1)/2$, and $v_i:=r+2i-1$, for $i\in[t]$. For any $v\in\sv$, we set 
$$A_v:=\{\sigma\in\hck^{n-2}\,|\,\sigma(v)=[n-1]\},\,\, B_v:=\{\sigma\in\hck^{n-3}\,|\,\sigma(v-1)\cup\sigma(v+1)=[n-1]\}.$$

For $S\subseteq\hck^d$, set $q(S):=\oplus_{\sigma\in S}\{\zz\sigma\}\in C^d(\hck/\zz)$, where $\zz\sigma=\{\sigma,\gamma\sigma\}$. We see that $q(A_0)=\emptyset$, and $q(S\oplus T)=q(S)\oplus q(T)$, for any $S,T\subseteq\hck^d$. Furthermore, since $\tau\cap\gamma\tau=\emptyset$, for any $\tau\in\hck^{d+1}$, we have $q(\cob S)=\cob q(S)$, for any $S\subseteq\hck^d$.  

It is easy to describe a~cochain representing $w_1^{n-2}(\hck)$. Let $\iota:K_2\hookrightarrow\sv$ be given by $\iota(1)=r$, $\iota(2)=r+1$, where $V(K_2)=[2]$. This induces an algebra homomorphism $\varphi:H^*(\thom(K_2,K_n)/\zz;\zz)\ra H^*(\hck/\zz;\zz)$. It is well-known that $\thom(K_2,K_n)/\zz\cong\rp^{n-2}$. Let $\tau\in\thom(K_2,K_n)^{n-2}$ be given by $\tau(1)=[n-1]$, $\tau(2)=\{n\}$. Since the dual of any cell generates $H^{n-2}(\rp^{n-2};\zz)$, we have $w_1^{n-2}(\thom(K_2,K_n))=[\{\zz\tau\}]$. By functoriality of $w_1$ we get $w_1^{n-2}(\hck)=[\varphi(\{\zz\tau\})]$. Comparing this to our notations we derive $w_1^{n-2}(\hck)=[q(A_r)]$.

\section{Proof of the Babson-Kozlov Conjecture}

\begin{lm} \label{lm:main}
We have $\cob B_v=A_{v-1}\oplus A_{v+1}$, for any $v\in\sv$.
\end{lm}
\pr The cells in $\cob B_v$ are obtained by taking a~cell $\sigma\in B_v$ and adding $x$ to $\sigma(w)$, for some $x\in[n]$, $w\in\sv$. When $w\neq v\pm 1$, we get a~cell $\tau$ which appears in $\cob B_v$ twice: in $\cob\sigma_1$ and in $\cob\sigma_2$, where $\sigma_1$, $\sigma_2$ are obtained from $\tau$ by deleting one of the elements from $\tau(w)$. When $w=v\pm 1$, we also get a~cell $\tau$ which appears in $\cob B_v$ twice: in $\cob\sigma_1$ and in $\cob\sigma_2$, where $\sigma_1$, $\sigma_2$ are obtained from $\tau$ by deleting $\{x\}=\tau(v-1)\cap\tau(v+1)$ either from $\tau(v-1)$ or from $\tau(v+1)$; unless $|\tau(v-1)|=1$ or $|\tau(v+1)|=1$. The latter cells appear once and yield $A_{v-1}\oplus A_{v+1}$.
\qed

\vspace{5pt}

\nin {\bf Proof of Theorem~\ref{thm:main}.}
Set $K:=\oplus_{i=1}^t q(B_{v_i})$, then 
$\cob K=\oplus_{i=1}^t \cob q(B_{v_i})=\oplus_{i=1}^t q(\cob B_{v_i})=\oplus_{i=1}^t(q(A_{v_i-1})\oplus q(A_{v_i+1}))=q(A_r)\oplus q(A_0)=q(A_r)$, hence $w_1^{n-2}(\hck)=[q(A_r)]=[\cob K]=0$.
\qed

\vspace{5pt}

We remark that the Babson-Kozlov Conjecture implies the Lov\'asz Conjecture, and that the latter was originally settled by Eric Babson and the author in~\cite{BK03a,BK03c}.



\end{document}